\documentclass[draft]{article}
\usepackage{amsmath,amsfonts,amssymb}
\usepackage[mathscr]{euscript}
\DeclareSymbolFont{bbold}{U}{bbold}{m}{n}
\DeclareSymbolFontAlphabet{\mathbbold}{bbold}

\newtheorem{theorem}{Theorem}[section]

\newtheorem{corollary}[theorem]{Corollary}
\newtheorem{lemma}[theorem]{Lemma}
\newtheorem{example}[theorem]{Example}
                 
\newtheorem{definition}[theorem]{Definition}

\newtheorem{remark}[theorem]{Remark}
\newtheorem{proposition}[theorem]{Proposition}

\def\mf{\mathfrak}\def\mc{\mathcal}\def\mb{\mathbb}\def\ms{\mathscr}
\def\R{\mb{R}}\def\N{\mb{N}}
\def\Q{\mb{Q}}\def\I{\mb{I}}
\def\S{\mb{S}}\def\T{\mb{T}}
\def\I{\mb{I}}\def\A{\mb{A}}\def\B{\mb{B}}
\def\s{\sigma}
\def\n{\wedge}
\def\ds{\displaystyle}

\def\P{\mathbb P}
\def\orth{\operatorname{Orth}}
\def\qed{\phantom{em}\hfill$\Box$}
\def\F{\mb{F}}
\def\L{\mb{L}}
\def\scr{\mathscr}
\numberwithin{equation}{section}
\def\pred{\operatorname{pred}}

\def\M{\mc{M}}
\def\dom{\operatorname{dom}}
\def\<{\langle}\def\>{\rangle}

\def\lr#1{\langle {#1}\rangle}
\def\wt{\widetilde}

\begin{document}
\title{Stopped processes and Doob's optional sampling theorem}
\author{Jacobus J. Grobler \\ Research Unit for Business Mathematics and Informatics\\
North-West University (Potchefstroom Campus),\\
Potchefstroom 2520,\\
South Africa\\
email: jacjgrobler@gmail.com\\ \\
Christopher Michael Schwanke \\
Department of Mathematics\\
Lyon College\\
Batesville, AR 72501, USA\\
and\\
Research Unit for Business Mathematics and Informatics\\
North-West University (Potchefstroom Campus),\\
Potchefstroom 2520,\\
South Africa\\
email: cmschwanke26@gmail.com
}
\maketitle
\begin{abstract} Using the spectral measure $\mu_\S$ of the stopping time $\S,$ we define the stopping element $X_\S$ as a Daniell integral $\int X_t\,d\mu_\S$ for an adapted stochastic process $(X_t)_{t\in J}$ that is a Daniell summable vector-valued function. This is an extension of the definition previously given for right-order-continuous sub-martingales with the Doob-Meyer decomposition property. The more general definition of $X_\S$ necessitates a new proof of Doob's optional sampling theorem, because the definition given earlier for sub-martingales implicitly used Doob's theorem applied to martingales. We provide such a proof, thus removing the heretofore necessary assumption of the Doob-Meyer decomposition property in the result.

Another advancement presented in this paper is our use of unbounded order convergence, which properly characterizes the notion of almost everywhere convergence found in the classical theory. Using order projections in place of the traditional indicator functions, we also generalize the notion of uniformly integrable sequences. In an essential ingredient to our main theorem mentioned above, we prove that uniformly integrable sequences that converge with respect to unbounded order convergence also converge to the same element in $\mc{L}^1$.
\end{abstract}


%%%%%%%%%%%%%%%%%%%%%Section 1 Introduction%%%%%%%%%%%%%%%%%%%%%%%%
\section{Introduction} In this paper we study the stopped process of a stochastic process in Riesz spaces. The notion of a stopped process is fundamental to the study of stochastic processes, since it is often used to extend results that are valid for bounded processes to hold also for unbounded processes. In~\cite{G2} we defined stopped processes for a class of submartingales and we expressed the need to get a definition applicable to more general processes. In this paper we introduce a much more general definition of stopped processes using the Daniell integral. The definition applies to every Daniell integrable process with reference to a certain spectral measure; this class of processes includes the important class of right-continuous processes.

\medskip
Considering the classical case, let $(\Omega,\mf{F},P)$ be a probability space, and let $(X_t=X(t,\omega))_{t\in J,\omega\in\Omega}$ be a stochastic process in the $L^1(\Omega,\mf{F},P)$ adapted to the filtration $(\mf{F}_t)$ of sub-$\sigma$-algebras of $\mf{F}.$ If the real valued non-negative stochastic variable $\S(\omega)$ is a stopping time for the filtration, then the stopped process is the process (see~\cite[Proposition 2.18]{KS}) 
$$
(X_{t\n \S})_{t\in J}=(X(t\n \S(\omega),\omega))_{t\in J,\, \omega\in\Omega}.
$$  
The paths of this process are equal to $X_t(\omega)$  up to time $\S(\omega),$ and from then on they remain constant with value $X_\S(\omega)=X(\S(\omega),\omega).$ 
\medskip 

The difficulty encountered in the abstract case is to define, what we shall call the {\em stopping element, $X_\S,$} needed in the definition of the stopped process $X_{t\wedge\S}.$ We note that $X_\S$ can be interpreted as an element of a vector-valued functional calculus on $\mf{E}$ induced by the vector function $t\mapsto X_t,$ in the same way as, in the case of a real-valued function $t\mapsto f(t)$ the element $f(X),$  $X\in\mf{E},$ is an element of the functional calculus on $\mf{E}$ induced by $f.$ The latter element can be obtained as a limit of simple elements of the form 
$\sum_{i=1}^n f(t_i)(E_{i+1}-E_i),$ with $E_i$ elements of the spectral system of $X$ with reference to a weak order unit $E$ (by taking $f(t)=t,$ the reader will recognize this as Freudenthal's spectral theorem). The element $f(X)$ can then be interpreted as an integral 
$$
\int_\R f(t)\,d\mu_E(t),
$$
with $\mu_E$ the spectral measure that is the extension to the Borel algebra in 
$\R$ of the (vector-) measure defined on left-open right-closed intervals $(a,b]$ by $\mu_E[(a,b]]=E_b-E_a$ (see~\cite[Sections IV.10, XI.5-7]{Vu}). Our approach will be to define a similar functional calculus for vector-valued functions, and we do this by employing the vector-valued Daniell integral as defined in~\cite{G12}. 

\medskip
Having a more general definition of $X_\S$ implies that a new proof of Doob's optional sampling theorem is needed. The reason is that a special case of  Doob's theorem (for martingales) is implicitly contained in the definition of $X_\S$ as given in \cite{G2}. Having obtained such a proof also proves that the definition given in \cite{G2} is a valid definition for the case considered there.

\medskip
A novelty in this paper is that we do not use order convergence in the definition of a continuous stochastic process, but unbounded order (uo-) convergence. This gives us a better model of the classic case, where convergence of a stochastic process is defined to mean that for almost every $\omega$ the paths  $X_t(\omega)$ of the process are continuous functions of $t.$ It is known that uo-convergence in function spaces is the correct notion to describe almost everywhere convergence. The fact, from integration theory, stating that a sequence that is pointwise almost everywhere convergent and uniformly integrable, is convergent in $L^1,$ is also generalized. This generalization is needed in the proof of Doob's optional sampling theorem. 

\medskip
We finally remark that we use, following \cite{G12, Pr}, the Daniell integral for vector-valued functions in our work. It turns out that Daniell's integral fits in perfectly in the Riesz space setting in which we describe stochastic processes. In fact, Daniell's original 1918 paper~\cite{PJD} was the first paper using, what we call today, Riesz space theory.  

%%%%%%%%%%%%%%%%%%Section 2 Preliminaries%%%%%%%%%%%%%%%%%%%%%%%%%%%%
\section{Preliminaries} Let $\mf{E}$ be a Dedekind complete, perfect Riesz space with weak order unit $E.$ We assume $\mf{E}$ to be separated by its order continuous dual $\mf{E}^\sim_{00}.$  For the theory of Riesz spaces (vector lattices) we refer the reader to the following standard texts~\cite{AB2,LZ,MN, Sch, Z1, Z2}. For results on topological vector lattices the standard references are~\cite{AB1,F}. We denote the {\it universal completion} of 
$\mathfrak E,$ which is an $f$-algebra that contains $\mathfrak E$ as an order dense ideal, by $\mathfrak{E}^u$ (the fact that it is an ideal follows from~\cite[Lemma 7.23.15 and Definition 7.23.19]{AB1}). Its multiplication is an extension of the multiplication defined on the principal ideal $\mf{E}_E$ and $E$ is the algebraic unit and a weak order unit for $\mathfrak E^u$ (see \cite{Z1}). The set of order bounded band preserving operators, called orthomorphisms, is denoted by $\operatorname{Orth}(\mf{E}).$  We refer to \cite{Donner, G9} for the definition and properties of the \emph{sup-completion} $\mathfrak{E}^s$ of a Dedekind complete Riesz space $\mathfrak{E}.$ It is a unique Dedekind complete ordered cone that contains $\mathfrak E$ as a sub-cone of its group of invertible elements and its most important property is that it has a largest element. Being Dedekind complete this implies that every subset of $\mf{E}^s$ has a supremum in $\mf{E}^s.$ Also, for every $C\in\mf{E}^s,$ we have $C=\sup\{X\in\mf{E}: X\le C\}$ and $\mf{E}$ is a solid subset of $\mf{E}^s.$

\medskip
A {\em  conditional expectation} $\mb{F}$ defined on $\mf{E}$ is a strictly positive order continuous linear projection with range a  Dedekind complete Riesz subspace $\mf{F}$ of $\mf{E}.$ It has the property that it maps weak order units onto weak order units. It may be assumed, as we will do, that $\mb{F}E=E$ for the weak order unit $E.$ The space $\mf{E}$ is called  {$\mb{F}$-universally complete} (respectively,  {$\mb{F}$-universally complete in $\mf{E}^u$})  if, whenever $X_\alpha\uparrow$ in $\mf{E}$ and $\mb{F}(X_\alpha)$ is bounded in $\mf{E}$ (respectively in $\mf{E}^u$), then $X_\alpha\uparrow X$ for some $X\in\mf{E}.$ If $\mf{E}$ is $\F$-universally complete in $\mf{E}^u,$ then it is $\F$-universally complete. 

\medskip

\textit{We shall henceforth tacitly assume that $\mf{E}$ is $\mb{F}$-universally complete in $\mf{E}^u.$}

\medskip For an order closed subspace $\mf{F}$ of $\mf{E},$ we shall denote the set of all order projections in $\mf{E}$ by $\mf{P}$ and its subset of all order projections mapping $\mf{F}$ into itself, by $\mf{P}_{\mf{F}}.$ This set can be identified with the set of all order projections of the vector lattice $\mf{F}$ (see~\cite{G1}).   

\medskip
B.A. Watson~\cite{W} proved that if $\mf{G}$ is an order closed Riesz subspace of $\mf{E}$ with $\mf{F}\subset\mf{G},$ then there exists a unique conditional expectation $\mb{F}_\mf{G}$ on $\mf{E}$ with range $\mf{G}$ and $\mb{F}\mb{F}_\mf{G}=\mb{F}_\mf{G}\mb{F}=\mb{F}$ (see~\cite{G2,W}). We shall also use the fact (see~\cite[Theorem 3.3 and Proposition 3.4]{G2}) that $Z=\mb{F}_\mf{G}(X),$ if and only if we have 
$$
\F(\P Z)=\F(\P X) \mbox{ holds for every projection $\P\in\mf{P}_{\mf{G}}$ }.
$$
 
The conditional expectation $\mb{F}$ may be extended to the sup-completion in the following way: For every $X\in\mf{E}^s,$ define $\mathbb F X$ by $\sup_{\alpha}\mathbb F X_\alpha\in\mathfrak{E}^s$ for any upward directed net  $X_\alpha\uparrow X$, $X_\alpha\in\mathfrak{E}.$ It is well defined (see~\cite{G8}). 
We define $\dom^+ \F:=\{0\le X\in\mf{E}^s:\ \F(X)\in\mf{E}^u\}.$ Then $\dom^+\F\subset\mf{E}^u$ (see~\cite[Proposition 2.1]{G6}) and we define $\dom\F=\dom^+\F-\dom^+\F.$ Since we're assuming $\mf{E}$ is $\F$-universally complete in $\mf{E}^u,$ we have $\dom\,\F=\mf{E}.$  

If $XY\in \rm{dom}\,{\mathbb F}$ (with the multiplication taken in the $f$-algebra $\mf{E}^u$), where $Y\in \mathfrak E$ and $X\in \mathfrak F=\mc{R}(\F),$ we have that $\mathbb F(XY)= X \mathbb F(Y)$. This fundamental fact is referred to as the \emph{averaging property} of $\mb F$ (see~\cite{G1}). 

Let $\Phi$ be the set of all $\phi\in\mf{E}^\sim_{00}$ satisfying $|\phi|(E)=1$ and extend $|\phi|$ to $\mf{E}^s$ by continuity.  Define $\mathscr{P}$ to be the set of all Riesz seminorms defined on $\mf{E}$ by $p_{\phi}(X):=|\phi|(\F(|X|)$ where $\phi\in\Phi.$ We define the space $\mathscr{L}^1:=(\mf{E},\s(\mathscr{P}))$ and have that $\mc{L}^1:=\{X\in\mf{E}^u: p_\phi(X)<\infty\mbox{ for all } \phi\in\Phi\},$ equipped with the locally solid topology $\s(\mathscr{L}^1,\mathscr{P})$ (for the proof see~\cite{G9}).  

We next define the space $\mc{L}^2$ to consist of all $X\in\mc{L}^1$ satisfying $|X|^2\in\mc{L}^1,$ where the product is taken in the $f$-algebra $\mf{E}^u.$ Thus, $\mc{L}^2:=\{X\in\mf{E}^u: |\phi|(\F(|X|^2))<\infty\mbox{ for all } \phi\in\Phi\}.$
For $X\in\mc{L}^2$ we define the Riesz seminorm 
$q_{\phi}(X):=(|\phi|(\F(|X|^2))^{1/2},$ and we denote the set of all these seminorms by $\mathscr{Q},$ and we equip $\mc{L}^2$ with the weak topology $\s(\mathscr{Q}).$ 

The spaces $\mathscr{L}^1$ and  $\mathscr{L}^2$ are topologically complete (see~\cite{G7} and \cite{G9} and note that this may not be true without the assumption that  $\mf{E}$ is $\F$-universally complete in $\mf{E}^u$).

\medskip
 A {\em filtration} on $\mf{E}$ is a set $(\F_t)_{t\in J}$ of conditional expectations satisfying $\F_s=\F_s\F_t$ for all $s<t.$ We denote the range of $\F_t$ by $\mf{F}_t.$ A {\em  stochastic process} in $\mf{E}$ is a function $t\mapsto X_t\in\mf{E},$ for $t\in J,$ with $J\subset\R^+$ an interval.  The stochastic process $(X_t)_{t\in J}$ is {\em adapted to the filtration} if $X_t\in\mf{F}_t$ for all $t\in J.$ 
 
 We shall write $\mf{P}_t$ to denote the set of all order projections that maps $\mf{F}_t$ into itself and we recall that $\F_t\P=\P\F_t$ holds for all $\P\in\mf{P}_t.$ The projections in $\mf{P}_t$ are the {\em events} up to time $t$ and $\mf{P}_t$ is a complete Boolean algebra.

Let $(\F_t)_{t\in J}$ be a filtration on $\mf{E}$. We recall that 
$$
\mf{F}_{t+}:=\bigcap_{s>t}\mf{F}_s,
$$ 
and the filtration is called right continuous if $\mf{F}_{t+}=\mf{F}_t$ for all $t\in J.$
Since we're assuming that $\mf{E}$ is $\F$-universally complete in $\mf{E}^u$, there exists a unique conditional expectation $\F_{t+}$ from $\mf{E}$ onto $\mf{F}_{t+}$ satisfying $\F\F_{t+}=\F_{t+}\F=\F$ (see~\cite[Proposition 3.8]{G2}). The set of order projections in the space $\mf{F}_{t+}$ will be denoted by $\mf{P}_{t+}.$

 If $(X_t)$ is a stochastic process adapted to $({\mathbb F}_t, \mathfrak F_t)$, we call $(X_t, \mathfrak F_t,\F_t)$ a \emph{supermartingale} (respectively \emph{submartingale}) if ${\mathbb F}_t(X_s)\leq X_t$ (respectively ${\mathbb F}_t(X_s)\geq X_t$) for all $t\leq s$. If the process is both a sub- and a supermartingale, it is called a \emph{martingale}. 
A stochastic process $(X_t)$ is said to be uo-convergent to $X\in\mf{E}$ as $t$ tends to $s,$ if 
$$
o-\lim_{t\to s}|X_t-X|\wedge Z=0
$$ 
for every positive $Z\in\mf{E}. $

We call a stochastic process uo-continuous in a point $s$ if 
 $$
 \operatorname{uo}-\lim_{t\to s}X_t=X_s.
 $$
In function spaces uo-convergence corresponds to pointwise almost everywhere convergence. Hence, the use of uo-convergence to define continuity in the abstract case, yields a direct generalization of the notion of path-wise continuity. The definitions of right-uo-continuity and left-uo-continuity in a point $s$ use uo-convergence from the right or from the left respectively.

\medskip
The band generated by $(tE-X)^+$ in $\mf{E}$ is denoted by 
$\mf{B}_{(tE>X)}$ and the projection of $\mf{E}$ onto this band by $\P_{(tE>X)}.$  The component of $E$ in $\mf{B}_{(tE>X)}$ is denoted by $E_t^\ell,$ i.e., $E_t^\ell=\P_{(tE>X)}E.$ The system
$(E_t^\ell)_{t\in J}$ is an increasing left-continuous system, called the {\em left-continuous spectral system of $X.$} Also, if $\overline{E}^r_t$ is the component of $E$ in the band generated by $(X-tE)^+$ and $E^r_t:=E-\overline{E}^r_t,$ the system $(E^r_t)$ is an increasing right-continuous system of components of $E,$ called the {\em right-continuous spectral system} of  $X$ (see~\cite{LZ,G2}).
The next definition was given in~\cite{G2}.

\begin{definition} \rm A {\em stopping time}  for the filtration 
$(\F_t,\mf{F}_t)_{t\in J}$ is an orthomorphism $\mb{S}\in\orth(\mf{E})$ such that its right continuous spectral system $(\mb{S}^r_t)$ of projections satisfies $\mb{S}^r_t\in\mf{P}_t.$ It is called an {\em optional time} if the condition holds for its left-continuous system $(\S_t^\ell)_{t\in J}.$ 
\end{definition}

We recall the fact that a stopping time is also an optional time and that the two concepts coincide for right-continuous filtrations. We shall use the following notation: $C^\ell_t:=\S_t^\ell E$ and $C^r_t:=\S^r_tE.$ The processes $C^\ell_\S:=(C^\ell_t)_{t\in J}$ and $C^r_\S:=(C^r_t)_{t\in J}$ are processes of components of $E.$ We have the following reformulation of the definition:

\begin{proposition} The orthomophism $\S$ is a stopping time for the filtration $(\F_t,\mf{F}_t)$ if and only if the stochastic process $C^r_\S$  is adapted to the filtration. Similarly, $\S$ is an optional time if and only if the stochastic process $C^\ell_\S$ is adapted to the filtration.  
\end{proposition}

\medskip
We recall (see~\cite{G2}) that the set of events (order projections) determined prior to the stopping time $\S$ is defined to be the family of projections
$$
\mf{P}_\S:=\{\P\in\mf{P}\,:\, \P\S^r_t\F_t=\F_t\P\S^r_t\mbox{ for all $t$}\}.
$$
$\mf{P}_\S$ is a complete Boolean sub-algebra of $\mf{P}$ and $\P\in\mf{P}_\S$ if and only if $\P\S^r_t\in\mf{P}_t$ for every $t\in J.$ The set  
$$
\mf{C}_\S:=\{\P E\,:\, \P\in\mf{P}_\S\},
$$
is a Boolean algebra of components of $E$ and we denote  by 
$\mf{F}_\S$ the order closed Riesz subspace of $\mf{E}$ generated by $\mf{C}_\S.$
By~\cite[Proposition 5.4]{G2}, there exists a unique conditional expectation $\F_\S$ that maps $\mf{E}$ onto $\mf{F}_\S$ with the property that $\F=\F\F_\S=\F_\S\F.$ 

Similarly, if $\S$ is an optional time for the filtration $(\F_t,\mf{F}_t),$ we find in~\cite{G2} that the Boolean algebra 
$\mf{P}_{\S+}$ of events determined immediately after $\S$ is the Boolean algebra of projections given by
\[
\mf{P}_{\S+}:=\{\P\in\mf{P}\, :\, \P\S_t^r\F_{t+} = \F_{t+}\P\S_t^r\mbox{ for all }t\}.
\] 
This is again a complete Boolean algebra of projections and 
$$
\mf{C}_{\S+}:=\{\P E\,:\, \P\in\mf{P}_{\S+}\}
$$
is a complete Boolean algebra of components of $E.$ We define the space $\mf{F}_{\S+}$ to be the Dedekind complete Riesz space generated by $\mf{C}_{\S+}.$ Since the space contains $\mf{F}_0,$  there exists a unique conditional expectation, denoted by 
$\F_{\S+},$ that maps $\mf{E}$ onto $\mf{F}_{\S+}$ with the property that $\F=\F\F_{\S+}=\F_{\S+}\F.$  

Our aim is now to define the stopping element $X_\S$ for a stopping time $\S,$ and having done that, the process $(X_{t\wedge\S})$ will be called the {\em stopped process.}

%%%%%%%%%%%%%%%%%Section 3 The Riesz space of simple processes%%%%%%%%%%%%%%
\section{Definition of the stopping element $X_\S$} 

Let $J=[a,b]$ and consider the optional time $\S$ for the filtration 
$(\F_t,\mf{F}_t)_{t\in J}$ defined on $\mf{E}$ with spectral interval contained in $J.$ Its left continuous spectrum of band projections $(\S^\ell_t)_{t\in J}$ is then adapted to the filtration meaning that $\S^\ell_t$ is a band projection in the Dedekind complete space $\mf{F}_t.$ Therefore,  the component $C^\ell_t$ of $E$ is an element of $\mf{F}_t.$ We define a vector measure $\mu_\S$ on the intervals $[s,t)$ by defining 
$$
\mu_\S[t_{i-1},t_i):=C^\ell_{t_i}-C^\ell_{t_{i-1}}.
$$
We refer the reader to~\cite[Section XI.5]{Vu} for a proof that this defines a 
$\sigma$-additive measure on the algebra of all finite unions of disjoint 
sub-intervals of the form $[s,t)$ in $J$ (and can be extended to $\sigma$-additive measure on the $\sigma$-algebra of Borel subsets of $J$).                       

\medskip
Next, let $\pi=\{a=t_0<t_1<\cdots<t_n=b\}$ be a partition of $J.$ We define $\L$ to be the Riesz space of all right-continuous simple processes of the form 
$$
X^\pi_t:=\sum_{t_i\in\pi}X_{t_i}\chi_{[t_{i-1}, t_i)},\ \ X_{t_i}\in\mf{F}_{t_i},
$$
where $\pi$ varies over all partitions of $J$ and $\chi_S$ is the indicator function of the set $S.$ 

We note that the process $(X^\pi_t)_{t\in J}$ is not in general adapted to the filtration $(\F_t,\mf{F}_t),$ because $X^\pi_{t_{i-1}}=X_{t_i}$ and $X_{t_i}$ need not be an element of $\mf{F}_{t_{i-1}}.$ We have that $X^\pi_{t}\in\mf{F}_{t_i}$ for all $t_{i-1}\le t<t_i.$ 

\medskip
The next proposition holds. 
\begin{proposition}
With pointwise ordering, $\L$ is a Riesz subspace of the \linebreak Dedekind complete Riesz space 
$\mf{E}^J$ of all $\mf{E}$-valued functions defined on the interval $J.$ 
\end{proposition}

\medskip
Let 
$$
I_\S(X_t^\pi):=\sum_{t_i\in\pi}X_{t_i}(C^\ell_{t_i}-C^\ell_{t_{i-1}})
=\sum_{t_i\in\pi}X_{t_i}\mu_\S[t_{i-1},t_i),
$$
where  the product  of $X\in\mf{E}$ and the component $C=\P E,$ $\P\in\mf{P},$ is defined to be $\P X.$ We can therefore also write 
$$
I_\S(X_t^\pi):=\sum_{t_i\in\pi}(\S^\ell_{t_i}-\S^\ell_{t_{i-1}})X_{t_i}.
$$
\begin{proposition} The operator $I_\S:\L\to\mf{E}$ has the following properties.
\begin{enumerate}
\item[\rm(1)] $I_\S$ is positive and linear;
\item[\rm(2)] $I_\S$ is $\sigma$-order continuous, i.e., if $X_{n,t}\in \L$ and $X_{n,t}\downarrow_n 0$ for each $t\in J,$ then $I_\S(X_{n,t})\downarrow_n 0;$ 
\item[\rm(3)] $I_\S$ is a lattice homomorphism.
\end{enumerate}
The operator $I_\S$ is therefore a positive vector-valued Daniell integral defined on $\L.$
\end{proposition}
{\em Proof.} Property (1) needs no proof. To prove (2), let $X_{n,t}\in\L$ satisfy $X_{n,t}\downarrow 0$ for every $t\in [a,b].$
Let 
$$
X_{1,t}=\sum_{i=1}^N\xi_{i}\chi_{[t_{i-1},t_i)}(t), \xi_i\in\mf{F}_i,
$$
and let $\epsilon>0$ be arbitrary. For each $i,$ $1\le i\le N,$ let 
$$
\mf{B}_{1,i}:=\mf{B}_{(\xi_{i}>\epsilon E)}
=\mf{B}_{(X_{1,t}\chi_{[t_{i-1},t_{i})}(t)>\epsilon E)},
$$ 
i.e., the band generated by 
$(\xi_{i}-\epsilon E)^+$ in $\mf{E}.$ We also define 
$$
\mf{B}_{n,i}:=\mf{B}_{(X_{n,t}\chi_{[t_{i-1},t_{i})}(t)>\epsilon E)} \mbox{  for }n\ge 1,
$$ 
and we denote the band projection onto $\mf{B}_{n,i}$ by 
$\P_{n,i}$ for each $i$ and $n.$ Then, since $X_{n,t}\downarrow$ for all $t, $ we have that 
$$
\mf{B}_{n+1,i}\subset \mf{B}_{n,i}\subset \mf{B}_{1,i}
\mbox{ for }1\le i\le N.
$$
For $n\ge 1,$ let $B_{n,i}\subset J$ be defined by
$$
B_{n,i}:=\{t\in[t_{i-1},t_i)\,:\, \P_{n,i}X_{n,t}>0\}.
$$
The definition of the simple function $X_{n,t}$ implies that each $B_{n,i}$ is a finite union of left-closed right-open intervals, and note also that by the definition of $\P_{n,i},$ we have for every $t\in B_{n,i},$ $\P_{n,i}X_{n,t}>\epsilon E.$ Since, for every fixed $t\in J$ we have that  $X_{n,t}\downarrow 0,$ it follows that  $B_{n,i}\downarrow 0.$ The vector measure $\mu_\S$ is a $\sigma$-additive measure and it follows, since $\mu_\S(J)=C^\ell_b-C^\ell_a\in\mf{E},$ that $\mu_\S(B_{n,i})\downarrow 0.$ Then,
\begin{align*}
I_\S(X_{n,t}\chi_{[t_{i-1},t_{i})}(t))
&=I_\S[(\P_{n,i}X_{n,t}+\P^d_{n,i}X_{n,t})(\chi_{B_{n,i}}(t)+ \chi_{B^c_{n,i}}(t))]\\
&\le (\P_{n,i}\xi_{i}+\epsilon E)\mu_\S(B_{n,i})+\epsilon(\P_i-\P_{i-1})E\\
&\le (\xi_{i}+\epsilon E)\mu_\S(B_{n,i})+\epsilon (\P_i-\P_{i-1})E.
\end{align*}
It therefore follows, for each $\epsilon>0$ and $i,$ $1\le i\le N,$ that 
$$
o-\lim_{n\to\infty}I_\S(X_{n,t}\chi_{[t_{i-1},t_i)}(t))\le \epsilon(\P_i-\P_{i-1})E.
$$ 
Summing over $i,$ we get
$$
\inf_nI_\S(X_{n,t})\le \epsilon(\P_b-\P_{a})E.
$$
This holds for every $\epsilon>0$ and so $I_\S(X_{n,t})\downarrow 0.$

We now prove (3). Let $X=X_t^\pi$ and $Y=Y_t^\pi$ be two elements of $\L$ written with the same partition $\pi$ of $J.$ Then, with $\Delta\S^\ell_{t_i}:=(\S^\ell_{t_i}-\S^\ell_{t_{i-1}})$
\begin{multline*}
I_\S(X\vee Y)=\sum_{i=1}^n \Delta\S^\ell_{t_i}(X_{t_i}\vee Y_{t_i})
=\sum_{i=1}^n \Delta\S^\ell_{t_i}X_{t_i}\vee \Delta\S^\ell_{t_i}Y_{t_i}\\
=\bigvee_{i=1}^n\Delta\S^\ell_{t_i}X_{t_i}\vee \Delta\S^\ell_{t_i}Y_{t_i}
=\bigvee_{i=1}^n\Delta\S^\ell_{t_i}X_{t_i}\vee \bigvee_{i=1}^n\Delta\S^\ell_{t_i}Y_{t_i}
=I_\S(X)\vee I_\S(Y).
\end{multline*}
Thus $I_\S$ is a Riesz homomorphism.\qed  

\medskip
Applying the Daniell extension procedure to the primitive positive integral $I_\S,$ we obtain an integral defined on the Riesz space $\mc{L}_\S$ of all Daniell $\S$-summable vector-valued functions that has the special property that it is a Riesz homomorphism.

\begin{theorem}
An adapted left-continuous process $(X_t)$ that is bounded by an $\S$-summable vector valued function $X$ is $\S$-summable. In particular, if $|X_t|\le ME,$ $M\ge 0,$ then $X=(X_t)$ is $\S$-summable. 
\end{theorem}
{\em Proof.} Let 
$\pi_n=\{a=t^{(n)}_0<t^{(n)}_1<\ldots<t^{(n)}_{2^n}=b\}$ be a diadic partition of $[a,b].$ Define the element $X_n$ by
\begin{equation}
X_n(t):=\sum_{i=1}^{2^n} X_{t^{(n)}_{i-1}}\chi_{[t^{(n)}_{i-1},t^{(n)}_i)}(t).
\end{equation}
Then $X_n(t)$ belongs to $\L$ and we claim that $X_n(t)$ converges to $X_t$ in every point $t\in[a,b].$ 

Fix an element $t_0\in[a,b].$ Then, for each $n$ we have that $t_0\in [t^{(n)}_{i-1},t^{(n)}_i)\in\pi_n$ for a unique $i,$ $1\le i\le 2^n,$ and $X_n(t_0)=X_{t^{(n)}_{i-1}}.$ If, at some stage, $t_0$ is the left endpoint of an interval $[t^{(n')}_{i-1},t^{(n')}_i),$ then, for all finer partitions, it will remain the left endpoint of some interval in that partition and so $X_n(t_0)=X_{t_0}$ for all $n\ge n'.$ We may therefore assume that $t_0>t^{(n)}_{i-1}$ for all $n$ (here, by abusing the notation, $t^{(n)}_{i-1}$ will always denote the left endpoint of the unique interval of $\pi_n$ to which $t_0$ belongs; $i$ will therefore also depend on $n$).   Since $t^{(n)}_i-t^{(n)}_{i-1}<(b-a)2^{-n},$ we have that $t^{(n)}_{i-1}\uparrow t_0$ as $n$ tends to infinity, and by the left continuity of $(X_t)$ we have that $X_n(t_0)=X_{t^{(n)}_{i-1}}$ converges to $X_{t_0}$ in order as $n\to\infty.$

Since $(X_t)$ is bounded by an $\S$-summable function $X,$ the Lebesgue domination theorem for the Daniell integral implies that $(X_t)$ is summable and that, with convergence in order,  
\begin{equation*}
\lim_{n\to\infty}I_\S (X_n(t))=I_\S (X_t).\tag*{$\Box$}
\end{equation*}

\begin{definition}\label{X_S} For $X\in\mc{L}_\S$ we define 
$$
X_\S:=I_\S(X_t).
$$
\end{definition}

For the proof of the next result we need the following fact about unbounded order convergence.

\begin{proposition}\label{proposition 3.4} Let $\mf{E}$ be a Dedekind complete Riesz space with weak order unit $E.$ Then the following are equivalent
\begin{itemize}
\item[\rm(1)] The sequence $(X_n)\subset\mf{E}$ is uo-convergent to $X\in\mf{E}.$
\item[\rm(2)] $(|X_n-X|\wedge kE)$ is order convergent to $0$ in $\mf{E}$ for every $k\in\N.$
\item[\rm(3)] The sequence $(X_n)$ is order convergent to $X$ in $\mf{E}^u.$
\end{itemize}
\end{proposition}
{\em Proof.} The implication (1)$\implies$(2) is clear.

(2)$\implies$(3): For every $k\in\N,$ there exists a sequence $(V_n^{(k)})$ in $\mf{E}$ such that $V_n^{(k)}\downarrow_n 0$ and 
$$
|X-X_n|\wedge kE\le V^{(k)}_n.
$$
Let  
$$
\mf{B}^{(k)}:=\bigcap_{n\in\N}\mf{B}_{(kE>|X-X_n|)}
$$ 
and let $\P_k$ be the projection onto $\mf{B}^{(k)}.$ We note that $\mf{B}^{(k)}$ is an increasing sequence of bands and so also is the sequence of projections $\P_k.$ We also have that 
$$
\P_k(|X-X_n|)=\P_k(|X-X_n|\wedge kE)\le\P_k V_n^{(k)}\mbox{ for all }n.
$$
We now put $\Q_1:=\P_1,\ \Q_2:=(\P_2-\P_1),\ldots \Q_n:=(\P_n-\P_{n-1}),\ldots.$ Then $(\Q_n)$ is a sequence of disjoint projections.
We define 
$$
V_n:=\sup_{k\in\N}\Q_kV_n^{(k)}\mbox{ for all }n.
$$
This supremum exists in $\mf{E}^u$ for each $n\in \N.$ Now,  since $E$ is a weak order we have for every $n\in\N$ that 
$$
|X-X_n|\wedge kE\uparrow |X-X_n|,
$$
and so 
$$
|X-X_n|=\sup_{k\in\N}\Q_k(|X-X_n|)\le\sup_{k\in\N} \Q_kV^{(k)}_n=V_n\downarrow 0\mbox{ in }\mf{E}^u.
$$
Hence, $(X_n)$ is order convergent to $X$ in $\mf{E}^u.$ 

(3)$\implies$(1): By assumption there exists a sequence $(Z_n)$  in $\mf{E}^u$  such that $|X-X_n|\le Z_n\downarrow 0.$ Consider an arbitrary $U\in\mf{E}.$ Then 
$$
|X-X_n|\wedge U\le Z_n\wedge U\in\mf{E},
$$
since, $\mf{E},$ being Dedekind complete, is an ideal in $\mf{E}^u$ (see our remark in section 2). Moreover, since $Z_n\downarrow 0$ in $\mf{E}^u,$ it follows that $Z_n\wedge U\downarrow 0$ in $\mf{E}$ by the order denseness of $\mf{E}$ in $\mf{E}^u.$
Therefore (1) holds.\qed

\begin{proposition} Let $(X^n_t)$ be a sequence in $\mc{L}_\S$ that is uo-convergent to $X_t\in\mc{L}_\S$ in each point $t\in J.$ Then $(X^n_\S)$ is uo-convergent to $X_\S.$ 
\end{proposition}
{\em Proof.} The constant vector-valued functions $t\mapsto kE$ are in $\L$ and therefore Daniell integrable. This shows that the sequence $|X^n_t-X_t|\n kE,$ which is order convergent to $0$ in each point $t,$ is pointwise bounded by the integrable function $t\mapsto kE.$ Therefore Lebesgue's dominated convergence theorem for the Daniell integral implies that $I_\S(|X^n_t-X_t|\wedge kE)$ is order convergent to $0.$ But, since $I_\S$ is a Riesz homomorphism, 
$$
I_\S(|X^n_t-X_t|\wedge kE)=|I_\S (X^n)-I_\S(X)|\wedge I_\S(kE)=
|I_\S (X^n)-I_\S(X)|\wedge kE,
$$
which shows that $X^n_\S=I_\S(X_n)\overset{uo}{\to}I_\S(X)=X_\S.$\qed 

\section{Uniform Integrability}
In this section we generalize the notion of uniform integrability. There are several ways in which one can do this, due to the different modes of convergence we have. It seems that convergence in $\mc{L}^1$ is the right notion to use in our case. The role of the integral is played by a conditional expectation $\F$ that is defined on the Dedekind complete Riesz space $\mf{E}.$ Our assumptions on $\mf{E}$ are as we stated them in section~2.
We recall that the Riesz semi-norm $p_\phi$ is defined as 
$$
p_\phi(X):=|\phi|(\F(|X|)),\ \ \phi\in\Phi
$$
and the topology of $\mc{L}^1$ is the locally solid topology 
$\sigma(\mc{L}^1,\mc{P}).$

\begin{definition}\label{UniformInt}\rm .
The sequence $(X_n)$ in $\mf{E}$ is called \textit{$\mc{L}^1$-uniformly integrable} whenever we have that, for every $p_\phi\in\mc{P},$
\begin{equation}\label{eqUniformInt}
p_\phi(\P_{(|X_n|\ge \lambda E)}|X_n|)\to 0
\mbox{ as $0\le\lambda\uparrow\infty$ uniformly in $n.$}
\end{equation}
This means that for every $\epsilon>0$ and for every $p_\phi\in\mc{P}$  there exists some $\lambda_0$ (depending on $\epsilon$ and $p_\phi$) such that for all $\lambda\ge\lambda_0,$ we have 
we have that
$$
p_\phi(\P_{(|X_n|\ge \lambda E)}|X_n|) <\epsilon \mbox{ for all $n\in\N.$}
$$

A stronger notion, that can be called \textit{order-uniformly integrable} is to have 
$$
\sup_{n\in\N}\F(\P_{(|X_n|\ge \lambda E)}|X_n|)\downarrow 0 \mbox{ as $\lambda\uparrow\infty$}.
$$ 
\end{definition}
If $(X_n)$ is order-uniformly integrable, we have for every $n$ that 
$$
p_\phi(\P_{(|X_n|\ge \lambda E)}|X_n|)\le p_\phi(\sup_{n\in\N}\F(\P_{(|X_n|\ge \lambda E)}|X_n|))\downarrow 0.
$$
It follows that $p_\phi(\P_{(|X_n|\ge \lambda E)}|X_n|)\to 0$ uniformly in $n$ as $\lambda\uparrow\infty.$ Thus, order-uniform integrability of $(X_n)$ implies $\mc{L}^1$-uniform integrability of $(X_n).$

\medskip
We note that for each fixed $n$ we have that $\P_{(|X_n|>\lambda E)}\downarrow 0$ as $\lambda\uparrow\infty.$ Therefore, for each fixed $n,$ $\P_{(|X_n|>\lambda E)}|X_n|\downarrow $0 as $\lambda\uparrow\infty$ and since $\F$ is order continuous, also $\F(\P_{(|X_n|>\lambda E)}|X_n|)\downarrow 0$ as $\lambda\uparrow\infty$ for each fixed $n.$ Therefore, if $(X_n)$ has only a finite number of non-zero elements, it is clear (see~\cite[Theorem 16.1]{LZ}) that $(X_n)$ is order-uniformly integrable.

\medskip
If $(X_n)$ is a bounded sequence in $\mc{L}^2,$ i.e., if for every $q\in \mc{Q}$ there exists a constant $M_\phi$ such that $q_\phi(X_n)\le M_\phi$ for all $n\in \N,$ then, by the Cauchy-inequality,
$$
p_\phi(\P_{(|X_n|>\lambda E)}|X_n|)\le q_\phi(\P_{(|X_n|>\lambda E)}E)q_\phi(X_n)\le q_\phi(\P_{(|X_n|>\lambda E)}E)M_\phi.
$$
Therefore, if $q_\phi(\P_{(|X_n|>\lambda E)}E)\to 0$ uniformly in $n$ as $\lambda\uparrow\infty,$ then $(X_n)$ is $\mc{L}^1$-uniformly integrable. The next proposition can be compared to~\cite[Theorem I.2.1]{DU}. 

\begin{proposition}\label{4P1} Let $0\le X\in\mf{E},$ let 
$(\P_{t})_{0\le t<\infty}$ be projections in $\mf{P}.$ Then, given any $p_\phi\in\mc{P}$ and $\epsilon>0,$ there exists a $\delta>0$ such that $p_\phi(\P_tE)<\delta,$ implies that 
$p_\phi(\P_{t} X)<\epsilon.$ Thus, if $p_\phi(\P_{t} E)$ converges to $0,$ then $p_\phi(\P_{t} X)$ converges to $0.$
\end{proposition}
{\em Proof.} Assume that the proposition is false. Then there exists an element $\phi\in\Phi$ and some $\epsilon>0$ such that, for every $k,$ there exists a projection $\P_{t_k}$ satisfying 
$p_\phi(\P_{t_k}E)<2^{-k}$ and  $p_\phi(\P_{t_k}X)>\epsilon.$ Define
$$
\Q_k:=\P_{t_k}\vee \P_{t_{k+1}}\vee\ldots,
$$
Then $\Q_k\downarrow$ and 
\begin{equation}\label{eq4.2.1}
p_\phi(\Q_kX)\ge p_\phi(\P_{t_k}X)>\epsilon.
\end{equation}
But,
$$
p_\phi(\Q_kE)\le \sum_{j=k}^\infty p_\phi(\P_{t_j} E)\le 2^{1-k} \downarrow 0, \mbox{ as $k\to\infty$}.
$$
Since $\F$ is strictly positive, $p_\phi=|\phi|\F$ is strictly positive on the carrier band $C_{\phi}$ of $\phi.$ Therefore, $\Q_kE \downarrow 0$ on $C_\phi.$ But then $\Q_kX\downarrow 0$ on $C_\phi$ and by the order continuity of $p_\phi,$ it follows that 
$p_\phi(\Q_kX)\downarrow 0.$ This contradicts  
(\ref{eq4.2.1}). \qed

\bigskip

The next theorem is also a generalization of a well-known fact about uniform integrability.

\begin{theorem}\label{TH43} The sequence $(X_n)$ in $\mf{E}$ is uniformly integrable if and only if it satisfies the following conditions: 
\begin{enumerate}
\item[{\rm(1)}]  $(X_n)$ is a bounded set in 
$\mc{L}^1$
\item[{\rm(2)}] For every $p_\phi\in\mc{P},$\ \ $p_\phi(\P|X_n|)\to 0$ uniformly in $n$ as $p_\phi(\P E)\to 0,$ i.e., given 
$\epsilon>0$ and $p_\phi\in\mc{P},$ there exists a $\delta>0$ such that, if $p_\phi(\P E)\le\delta,$ then $p_\phi(\P|X_n|)<\epsilon$
for all $n\in\N.$
\end{enumerate}
\end{theorem}
{\em Proof.} Suppose that $(X_n)$ is a bounded set in $\mc{L}^1$ and that it is uniformly continuous, i.e., that $(X_n)$ satisfies condition (2). Then we have by Chebyshev's inequality, we have  
$$
\F(\P_{(|X_n|\ge tE)}E)\le \frac{1}{t}\F(|X_n|),
$$
which implies that 
$$
p_\phi(\P_{(|X_n|\ge tE)}E)\le\frac{1}{t}p_\phi(X_n)\le M_\phi,
$$
for a number $M_\phi\ge 0.$ By the boundedness of $(X_n)$ in $\mc{L}^1,$ it follows that $p_\phi(\P_{(|X_n|\ge tE)}E)\to 0$ uniformly in $n.$ It follows by (2) that 
$$
p_\phi(\P|X_n|)\to 0 \mbox{ uniformly in $n$}
$$
and so $(X_n)$ is uniformly integrable.

\medskip Conversely, if $(X_n)$ is uniformly integrable, we have for every $p_\phi\in\mc{P}$ that 
\begin{align}\label{equation4.3}
p_\phi(\P|X_n|)&=p_\phi(\P\P_{(|X_n|\ge tE)}|X_n|)
+p_\phi(\P\P_{(|X_n|< tE)}|X_n|) \nonumber\\
&\le p_\phi(\P_{(|X_n|\ge tE)}|X_n|) + tp_\phi(\P E).
\end{align}
By the uniform integrability, we can choose, for given $\epsilon>0,$ a number $t_0$ such that the first term is less that $\epsilon/2$ for all $n.$ We then have, for $\phi(\P E)<\epsilon/2t_0$ that $p_\phi(\P|X_n|)<\epsilon$ for all $n,$ thus proving that condition (2) holds.

\medskip\noindent Taking the projection $\P$ in (\ref{equation4.3})
equal to the identity $I,$ it follows that for large $t$ (depending on $p_\phi$) we have
$$
p_\phi(|X_n|)\le \epsilon+t=M_\phi<\infty.
$$
Since this holds for arbitrary $p_\phi\in\mc{P},$ the set $(X_n)$ is bounded in $\mc{L}^1.$\qed

\bigskip

\begin{corollary}\label{L:Xn+X}
If $(X_n)$ and $(Y_n)$ are uniformly integrable sequences, then $(X_n+Y_n)$ is also uniformly integrable. In particular, if $X\in\mf{E}$ then $(X_n+X)$ is uniformly integrable.
\end{corollary}
{\em Proof.} It is clear that if $(X_n)$ and $(Y_n)$ are bounded sequences in $\mc{L}^1$ then $(X_n+Y_n)$ is also a bounded sequence in $\mc{L}^1.$ Also, since they are uniformly integrable, they are uniformly continuous, i.e., condition (2) in Theorem 4.3 above holds for both of them. But then, for every $p_\phi\in\mc{P},$ we have that if $p_\phi(\P E)\to 0,$ then 
$$
p_\phi(\P|X_n+Y_n|)\le p_\phi(\P|X_n|)+p_\phi(\P |Y_n|)\to 0
$$
uniformly in $n.$ By Theorem~\ref{TH43}, this implies that $(X_n+Y_n)$ is uniformly integrable. \qed

\bigskip
Below we denote unbounded order convergence of a sequence $(X_n)$ to an element $X$ by $X_n\overset{uo}{\to}X.$

\begin{lemma}\label{L:zero}
If $X_n\overset{uo}{\to}0$ and $(X_n)$ is uniformly integrable, then $X_n\to 0$ in $\mc{L}^1.$
\end{lemma}

{\em Proof.} Suppose that $X_n\overset{uo}{\to}0$ and that $(X_n)$ is uniformly integrable. Let $\epsilon>0$ and $p_\phi\in\mc{P}$ be given. Then, it follows from the uniform integrability, that 
\begin{align*}
p_\phi(X_n)&=p_\phi(\P_{(|X_n|\geq\lambda E)}X_n)
+p_\phi(\P_{(|X_n|<\lambda E))}X_n)\\
&\le p_\phi(\P_{(|X_n|\geq\lambda E)}X_n)
+p_\phi(|X_n|\wedge \lambda E)\\
&<\epsilon/2 + p_\phi(|X_n|\wedge \lambda_0 E),
\end{align*}
for some $\lambda_0>0$ and for all $n\in\N.$ Since $X_n\overset{uo}{\to}0$ by assumption, and since $p_\phi$ is order continuous,  there exist some $N\in\N$ such that for all $n\ge N$ we have that the last term above is less that $\epsilon/2.$ Thus, $p_\phi(X_n)\to 0$ and this holds for every $p_\phi\in\mc{P}.$ \qed 

%However, using the well-known identity $\P_{g^+}(g)=g^+$ (see, %e.g. p. 215 in Introduction to Operator Theory in Riesz Spaces), %we have
%\begin{align*}
%\P(|X_n|<\lambda E)(\lambda E)-\P(|X_n|<\lambda E)(|X_n|)&=\P(|%X_n|<\lambda E)(\lambda E-|X_n|)\\
%&=(\lambda E-|X_n|)^+\\
%&\geq 0.
%\end{align*}
%It follows that $\P(|X_n|<\lambda E)(|X_n|)\leq\P(|X_n|<\lambda %E)(\lambda E)$ and therefore
%\begin{align*}
%&\P(|X_n|<\lambda E)(|X_n|)\\
%&=\P(|X_n|<\lambda E)(|X_n|)\wedge\P(|X_n|<\lambda E)(\lambda E)%\\
%&=\P(|X_n|<\lambda E)(|X_n|\wedge\lambda E).
%\end{align*}
%Hence we obtain
%\begin{align*}
%|\F(X_n)|&\leq Z_m+\F(\P(|X_n|<\lambda E)(|X_n|\wedge\lambda E))%\\
%&\leq Z_m+\F(|X_n|\wedge\lambda E).
%\end{align*}
%By assumption we have $|X_n|\wedge\lambda E\overset{o}{\to}0$. %It now follows from the order continuity of 
%$\F$ that $\F(X_n)%\overset{o}{\to}0$ as well.
%\qed

\begin{theorem}\label{T:FXn->FX}
If $X_n\overset{uo}{\to}X$ and $(X_n)$ is uniformly integrable, then $X_n\to X$ in $\mc{L}^1.$
\end{theorem}
{\em Proof.} Suppose that $X_n\overset{uo}{\to}X$ and that $(X_n)$ is uniformly integrable. For each $n\in\mb{N}$ define $C_n:=X_n-X$. Then $C_n\overset{uo}{\to}0$, and by Corollary~\ref{L:Xn+X}, we know that $(C_n)$ is uniformly integrable. Thus, by Lemma~\ref{L:zero}, it is true that $C_n\to 0$ in $\mc{L}^1.$ But this is equivalent to $X_n\to\, X$ in $\mc{L}^1.$
\qed

\bigskip
\noindent
\textbf{Conclusion} If the sequences $(X_{\S_n})$ and $(X_{\T_n})$ are uniformly integrable and uo-converge in $\mf{E}$ to $X_\S$ and $X_\T$, respectively, then it is easy to see that for any band projection $\P$ we have that $\P(X_{\S_n})\overset{uo}{\to}\P(X_\S)$ and  $\P(X_{\T_n})\overset{uo}{\to}\P(X_\T)$. It is also easy to see that $(\P(X_{\S_n}))$ and $(\P(X_{\T_n}))$ are also uniformly integrable. Therefore, by Theorem~\ref{T:FXn->FX} we have
$X_{\S_n}\to X_\S$ and $\P X_{\T_n}\to\P X_\T$ in $\mc{L}^1.$
This fact will be used in the proof of Doob's optional sampling theorem below.

\begin{definition}\rm (see~\cite[Problem 3.11]{KS}) Let $(\mf{F}_n)$ be a decreasing sequence of Dedekind complete Riesz subspaces of $\mf{E},$ i.e., 
$$
\mf{F}_{n+1}\subseteq \mf{F}_{n}\subseteq \mf{E},
$$
with $\mf{F}_n$ the range of a conditional expectation $\F_n:\mf{E}\to\mf{F}_n$ satisfying $\F_n\F_m=\F_m\F_n=\F_m$ if $m>n.$ 
The process $(X_n)$ with $X_n\in\mf{F}_n$ and $\F_{n+1}(X_n)\ge X_{n+1}$ is called a {\em backward submartingale}. 
\end{definition}
We note that $\mf{F}_\infty:=\bigcap_n\mf{F}_n$ is a Dedekind complete Riesz space that is contained in each of the spaces $\mf{F}_n$ and so there exists a conditional expectation $\F_\infty:\mf{E}\to\mf{F}_\infty$ with the property that for each $n$ we have $\F_\infty\F_n=\F_n\F_\infty=\F_\infty.$ Furthermore, applying $\F_\infty$ to both sides of the inequality in the definition, we find that for all $n,$ $\F_\infty(X_n)\ge\F_\infty(X_{n+1}),$ i.e., the sequence $(\F_\infty(X_n))$ is a decreasing sequence. It is also easy to show by induction that for all $n$ one has $\F_n(X_1)\ge X_n.$

\begin{example}\rm Let $(X_t,\F_t)_{t\in J}$  be a submartingale. With 
$J=[a,b],$ we have for any sequence of real numbers $t_n\downarrow a,$ that $(X_{t_n},\F_{t_n})_{n\in \N})$ is a backward submartingale. In this case, 
$\F_\infty=\F_a=\F.$
\end{example}

Since we work in the setting where we do not have an integral, but a fixed conditional expectation $\F,$ we shall assume for all backward submartingales considered that $\F_\infty=\F.$

\medskip
\begin{proposition}\label{4P2} Let $(X_n)$ be a backward submartingale with $\F_\infty=\F.$ If  the sequence  
$(\F(X_n))$ is bounded below, i.e., if 
$$
Y=\inf_{n\in\N}\F(X_n)\mbox{ exists in }\mf{E},
$$
then the sequence $(X_n)$ is uniformly integrable.
\end{proposition} 
{\em Proof.} %(The proof is analogous to that given in \cite[3.11, p.41]{KS} for %the classical case.) 
By Jensen's inequality, $(X_n^+, \mf{F}_n)$  is also a backward submartingale. Hence, for $\lambda>0,$ we find by the Chebyshev inequality that for each $n,$
$$
\lambda\F(\P_{(|X_n|>\lambda E)}E)\le\F(|X_n|)
=-\F(X_n)+2\F(X^+_n)\le -Y+2\F(X_1^+).
$$
It follows that
\begin{equation}\label{4E1}
\lim_{\lambda\to\infty}p_\phi(\P_{(|X_n|>\lambda E)}E)=0 \mbox{ uniformly in $n,$} 
\end{equation}
and therefore also 
\begin{equation}\label{4E2}
\lim_{\lambda\to\infty}p_\phi(\P_{(X_n^+>\lambda E)}E)=0 \mbox{ uniformly in $n.$}
\end{equation}
Using the backward submartingale property of $(X_n^+),$ we have 
\begin{multline}\label{4E3}
\F(\P_{(X_n^+>\lambda E)}X_n^+)
\le\F(\P_{(X_n^+>\lambda E)}\F_n X_1^+)\\
=\F\F_n(\P_{(X_n^+>\lambda E)}X_1^+)
=\F(\P_{(X_n^+>\lambda E)}X_1^+).
\end{multline}
Hence, we have for any $p_\phi\in\mc{P},$ that
\begin{equation}\label{E47}
p_\phi(\P_{(X_n^+>\lambda E)}X_n^+)\le p_\phi(\P_{(X_n^+>\lambda E)}X_1^+).
\end{equation}
We now apply Proposition \ref{4P1} to find for every $\epsilon>0,$ a $\delta>0$ such that, if
$p_\phi(\P_{(X_n^+>\lambda E)}E)<\delta,$ then $p_\phi(\P_{(X_n^+>\lambda E)}X_1^+)<\epsilon.$ From 
(\ref{4E2}), there exist some $\lambda_0$ such that, for $\lambda>\lambda_0,$ $p_\phi(\P_{(X_n^+>\lambda E)}E)<\delta,$ for all $n\in\N.$ It then follows from (\ref{E47}) that for all $\lambda>\lambda_0,$ we have $p_\phi(\P_{(X_n^+>\lambda E)}X_n^+)<\epsilon$ for all $n\in \N.$ This shows that the backwards submartingale $(X_n^+)$ is uniformly integrable.

\medskip
We next show that the sequence $(X_n^-)$ is also uniformly integrable. Note that $\P_{(X_n^->\lambda E)}=\P_{(X_n<-\lambda E)}$ and that for $m<n,$ we have $X_n\le\F_nX_m.$ Now,\begin{multline}\label{4E4}
0\ge\F(\P_{(X_n<-\lambda E)}X_n)=\F(X_n)-\F(\P_{(X_n\ge-\lambda E)}X_n)\\
\ge\F(X_n)-\F(\P_{(X_n\ge-\lambda E)}\F_nX_m)\\
\ge\F(X_n)-\F(\P_{(X_n\ge-\lambda E)}X_m)\\
=\F(X_n)-\F(X_m)+\F(\P_{(X_n<-\lambda E)}X_m).
\end{multline}  
Since the sequence $\F(X_n)\downarrow_n Y,$ $(X_n)$ is convergent in $\mc{L}^1$ and therefore a Cauchy sequence. For a given 
$\epsilon>0,$ we can choose $m=m(\epsilon)$ such that for all $n>m,$ we have 
$$
p_\phi(X_m-X_n)<\epsilon/2.
$$

Also, by Proposition~\ref{4P1}, there exists a $\delta>0,$ such that 
$p_\phi(\P_{(|X_n|>\lambda E)}E)<\delta,$ implies that 
$p_\phi(\P_{(|X_n|>\lambda E)}X_m)<\epsilon/2$ and using (\ref{4E1}), we can find a $\lambda_0$ such that for all $\lambda>\lambda_0,$ we have for all $n\in \N$ that 
$p_\phi(\P_{(|X_n|>\lambda E)}E)<\delta$ and therefore for all $n\in\N$ that $p_\phi(\P_{(|X_n|>\lambda E)}X_m)<\epsilon/2$ if $\lambda>\lambda_0.$ But, 
$\P_{(X_n^->\lambda E)}\le \P_{(|X_n|>\lambda E)}$ and therefore, there exists a $\lambda$ such that for all $\lambda>\lambda_0,$
$$
p_\phi(\P_{(X_n^->\lambda E)}X_m)<\epsilon/2 \mbox{ for all $n\in\N.$}
$$
We now use the inequality in (\ref{4E4}): For all $n>m(\epsilon)$ we have 
\begin{align}\label{4E5}
\F(\P_{(X_n^->\lambda E)}X_n^-)
&=|\F(\P_{(X_n^->\lambda E)}X_n)|\nonumber \\
&=-\F(\P_{(X_n<-\lambda E)}X_n) \\
&\le(\F(X_m)-\F(X_n))-\F(\P_{(X_n<-\lambda E)}X_m)
\end{align}
and so, for all $n>m$ we get
$$
p_\phi(\P_{(X_n^->\lambda E)}X_n^-)\le p_\phi(X_m-X_n)+
p_\phi(\P_{(X_n^->\lambda E)}X_m)<\epsilon/2+\epsilon/2=\epsilon
$$
for all $\lambda>\lambda_0.$

For $n=1,2,\ldots,m$ we have that 
$p_\phi(\P_{(X_n^->\lambda E)}X_n^-)\downarrow 0$ as $\lambda\to\infty$ so we can find $\lambda_n$ such that for $\lambda>\lambda_n,$ we have $p_\phi(\P_{(X_n^->\lambda E)}X_n^-)<\epsilon.$ If $\lambda>\max\{\lambda_0,\lambda_1,\ldots,\lambda_m\}$ we have that 
$$
p_\phi(\P_{(X_n^->\lambda E)}X_n^-)<\epsilon \mbox{ for all $n\in\N$}
$$
Thus, $(X_n^-)$ is uniformly integrable.
Our final result, that $(X_n)=(X_n^+-X_n^-$) is uniformly integrable, follows from Corollary \ref{L:Xn+X}.\qed

%%%%%%%%%%%%%%Section 5 The optional sampling theorem%%%%%%%%%%%%%%%%%%%%%
\section{The optional sampling theorem}
As remarked in the introduction, we have to prove the optional sampling theorem using Definition~\ref{X_S}. 
  
\begin{theorem} Let $(X_t)_{t\in J}$ be a right-uo-continuous submartingale and let $\S\le\T$ be two optional times of the filtration $(\F_t,\mf{F}_t).$ Then, if either 
\begin{enumerate}
\item $\T$ is bounded or
\item $(X_t)$ has a last element,
\end{enumerate}
we have
$$
\F_{\S+}X_\T\ge X_\S.
$$ 
If $\S$ and $\T$ are stopping times, one has
$$
\F_{\S}X_\T\ge X_\S.
$$
\end{theorem}
{\em Proof.} Let $\pi_n=\{a=t_0<t_1<\ldots<t_{2^n}=b\}$ be a diadic partition of $J=[a,b]$ and define the sequence $(\S_n)$ by putting
\begin{equation}\label{5E1}
\S_n=\sum_{i=1}^{2^n}t_i(\S^\ell_{t_i}-\S^\ell_{t_{i-1}})
=\sum_{i=1}^{2^n}t_i\Delta\S^\ell_i=\sum_{i=1}^{2^n}t_i\S^\ell_{t_i}(\S_{t_{i-1}}^{\ell})^d
\end{equation}
and similarly,
\begin{equation}\label{5E2}
\T_n=\sum_{j=1}^{2^n}t_j(\T^\ell_{t_j}-\T^\ell_{t_{j-1}})
=\sum_{j=1}^{2^n}t_j\Delta\T^\ell_j\sum_{i=1}^{2^n}t_j\T^\ell_{t_j}(\T_{t_{j-1}}^{\ell})^d.
\end{equation}
We now write them both as sums with respect to the partition $\{\Delta\S^\ell_i\Delta\T^l_j\}_{i,j=1}^n,$ i.e., we get
\begin{equation}\label{5E3}
\S_n=\sum_{i=1}^{2^n}\sum_{j=1}^{2^n}s_{ij}\Delta\T^\ell_j\Delta\S^\ell_i,\ s_{ij}=t_i\mbox{ and }  \T_n=\sum_{i=1}^{2^n}\sum_{j=1}^{2^n}t_{ij}\Delta\T^\ell_j\Delta\S^\ell_i,\ t_{ij}=t_j.
\end{equation}
Now, $\S\le\T$ implies that for each fixed $n,$ $\S_n\le\T_n$ and so $s_{ij}\le t_{ij};$ this implies that $t_i\Delta\T^\ell_j\Delta\S^\ell_i\le t_j\Delta\T^\ell_j\Delta\S^\ell_i$ for all $i,j$ such that $\Delta\T^\ell_j\Delta\S^\ell_i\ne 0.$

Each $\S_n$ and each $\T_n$ is a stopping time for the filtration and by Freudenthal's theorem, $\S_n\downarrow \S$ and $\T_n\downarrow\T.$ 

With these definitions for $\S_n$ and $\T_n,$ we have that 
\begin{equation}\label{5E4}
X_{\S_n}=\sum_{i=1}^{2^n}\sum_{j=1}^{2^n}\Delta\T^\ell_j\Delta\S^\ell_i X_{t_{i}}\mbox{ and }
X_{\T_n}=\sum_{i=1}^{2^n}\sum_{j=i}^{2^n}\Delta\T^\ell_j\Delta\S^\ell_i X_{t_{j}}.
\end{equation}
Next, we put 
\begin{equation}\label{5E5}
\F_{\S_n}:=\sum_{i=1}^{2^n}\sum_{j=1}^{2^n}\F_{t_{i}}\Delta\T^\ell_j\Delta\S^\ell_i. 
=\sum_{i=1}^{2^n}\F_{t_{i}}\Delta\S^\ell_i.
\end{equation}
It is readily checked that $\F_{\S_n}$ is a strictly positive, order continuous projection that maps $E$ onto $E,$ i.e., $\F_{\S_n}$ is a conditional expectation. Its range is the direct sum 
$$
\bigoplus_{i=1}^{2^n}\Delta\S^\ell_i\mf{F}_{t_i}
$$
of the bands $\Delta\S_i^\ell \mf{F}_{t_i}$ (since $\F_{t_i}$ and $\Delta\S_i^\ell $ commute), and the projections in $\mf{E}$ that belong to this space are exactly those projections in $\mf{E}$ that belong to $\mf{P}_t$ for all $t$ such that $\S_n \le tE.$ Therefore, the space $\mf{F}_{\S_n},$ which is by  definition the space generated by these projections, is equal to the space  
$\bigoplus_{i=1}^{2^n}\Delta\S^\ell_i\mf{F}_{t_i}.$ Moreover,
$$
\F\F_{\S_n}=\sum_{i=1}^{2^n}\F\F_{t_{i}}\Delta\S^\ell_i=\F\sum_{i=1}^{2^n}\Delta\S^\ell_i=\F,
$$
and similarly $\F_{\S_n}\F=\F.$ Therefore, $\F_{\S_n}$ is the unique conditional expectation with range $\mf{F}_{\S_n}$ satisfying these two conditions. 

For each fixed $i,$ consider the sum 
\begin{align}\label{equation5.6a}
\sum_{j=i}^{2^n}\F_{t_{i}}\Delta\T^\ell_j\Delta\S^\ell_i X_{t_j}
&=\F_{t_i}\Delta\T^\ell_{2^n}\Delta\S^\ell_i X_{t_{2^n}}+\ldots 
+\F_{t_i}\Delta\T^\ell_{t_i}\Delta\S^\ell_i X_{t_{i}}
\end{align}
and note that, in the first term, 
$$
\Delta\T^\ell_{2^n}\Delta\S^\ell_i=\Delta \S^\ell_i(\T^\ell_{t_{2^n}}(\T^{\ell}_{t_{2^n-1}})^d)=\Delta \S^\ell_i(\T^{\ell}_{t_{2^n-1}})^d\in\mf{P}_{t_{2^n-1}}.
$$
Therefore,
\begin{align*}
\F_{t_i}\Delta\T^\ell_{2^n}\Delta\S^\ell_i X_{t_{2^n}}&=\F_{t_i}\F_{t_{2^n-1}}
\Delta\T^\ell_{2^n}\Delta\S^\ell_i X_{t_{2^n}}\\
&=\F_{t_i}\Delta\T^\ell_{2^n}\Delta\S^\ell_i \F_{t_{2^n-1}}
X_{t_{2^n}}\\
&\ge\F_{t_i}\Delta\T^\ell_{2^n}\Delta\S^\ell_i X_{t_{2^n-1}}.
\end{align*}
Substituting this inequality in Equation (\ref{equation5.6a}), and repeating the process we finally arrive at 
\begin{equation}
\sum_{j=i}^{2^n}\F_{t_{i}}\Delta\T^\ell_j\Delta\S^\ell_i X_{t_j}
\ge\F_{t_{i}}\Delta\S^\ell_i(\sum_{j=i}^{2^n}\Delta\T^\ell_j) X_{t_i}
=\F_{t_{i}}\Delta\S^\ell_i(\sum_{j=i}^{2^n}\Delta\T^\ell_j) X_{t_i}
=\Delta\S^\ell_iX_{t_i}. 
\end{equation}
Thus, 
\begin{equation}\label{5E6}
\F_{\S_n}(X_{\T_n})=\sum_{i=1}^{2^n}\sum_{j=i}^{2^n}\F_{t_{i}}\Delta\T^\ell_j\Delta\S^\ell_i X_{t_j}
\ge \sum_{i=1}^{2^n}\Delta\S^\ell_i X_{t_i}=X_{\S_n}.
\end{equation} 
(This is Doob's optional sampling theorem for this special case.) 

For every $\P\in\mf{P}_{\S_n},$ we therefore have that
\begin{equation}\label{5E7}
\F(\P X_{\T_n})=\F\F_{\S_n}\P X_{\T_n}= \F\P(\F_{\S_n}X_{\T_n})\ge\F(\P X_{\S_n}).
\end{equation}
By~\cite[Proposition 5.15]{G2}, we have that 
$\displaystyle \mf{P}_{\S+}=\bigcap_{n=1}^\infty \mf{P}_{\S_n}.$
Therefore, by (\ref{5E7}), 
\begin{equation}\label{5E7b}
\F(\P X_{\T_n})\ge\F(\P X_{\S_n})\mbox{  holds for all $\P\in\mf{P}_{\S+}.$}
\end{equation} If $\S$ is a stopping time, it follows from~\cite[Proposition 5.9]{G2} that since $\S\le\S_n,$ we have $\mf{P}_{\S}\subset\mf{P}_{\S_n}$ and so this inequality holds in that case also for all $\mf{\P}\in\mf{P}_\S.$ 

\medskip
Applying the arguments above to $\S_n\le\S_{n+1}$ we get as in (\ref{5E6}) that  
$$
\F_{\S_n}(X_{\S_{n+1}})\ge X_{\S_n}\mbox{ for all $n$},
$$
which implies that $(\S_n,\F_{\S_n})$ is a backward submartingale and so
$(\F(X_{\S_n}))$ is a decreasing sequence and, using (\ref{5E4}), $\F(X_{\S_n})\ge\F(X_a)$ for all $n.$ Applying Proposition~\ref{4P2}, we have that the sequence $(X_{\S_n})$ is uniformly integrable. The same is true for the sequence $(X_{\T_n}).$

\medskip We now note that 
$$
X_{\S_n}=I_\S(X^n),\mbox{ and }X_{\T_n}=I_{\T}(X^n),
$$
with 
$$
X^n_t=\sum_{i=0}^{2^n}X_{t_i}\chi_{[t_{i-1},t_i)}(t).
$$
Our assumption that $X=(X_t)$ is uo-right-continuous, implies that in each point $t$ we have that $X^n_t$ is uo-convergent to $X_t.$ Now let $k\in\N$ and define the constant process $(kE_t=E).$ This process is Daniell integrable since 
$I_\S(kE)=kE\in\mf{E}.$ Consider the process $X^n\wedge kE=(X^n_t\wedge kE_t)=(X^n_t\wedge kE).$ Then $X^n\wedge kE$ converges in order pointwise to  $X\wedge kE.$ By Lebesgue's dominated convergence theorem, we get that 
$ I_\S(X^n\wedge kE) $ is order convergent to $I_\S(X\wedge kE).$ But, $I_\S$ is a Riesz homomorphism, so we get that 
$$
I_\S(X^n\wedge kE)=I_\S(X^n)\wedge I_\S(kE)=I_\S(X^n)\wedge kE.
$$ 
Therefore, by Proposition~(\ref{proposition 3.4}), $X_{\S_n}=I_\S(X^n)$ is uo-convergent to $I_\S(X)=X_{\S}$ and the same holds for $I_\T(X^n)$ and $I_\T(X).$ By our result for a uo-convergent sequence that is uniformly integrable we have that $p_\phi(\P X_{\S_n})$ converges to $p_\phi(\P X_\S)$ and also $p_\phi(\P X_{\T_n})$ converges to $p_\phi(\P X_\T),$ for every $\P\in\mf{P}_{\S+}$ and for every $p_\phi\in\mc{P}.$ Recalling that $p_\phi=|\phi|\F$ for $\phi\in\Phi$, this implies, using~(\ref{5E7b}), that  
$$
|\phi|\F(\P X_\T)\ge |\phi|\F(\P X_\S),\mbox{ for every $\P\in\mf{P}_{\S+}.$}
$$
But since $\mf{E}^\sim_{00}$ separates the points of $\mf{E},$ we get 
$$
\F(\P X_\T)\ge \F(\P X_\S),\mbox{ for every $\P\in\mf{P}_{\S+}$},
$$
and thus
$$
\F\F_{\S+}(\P X_\T)=\F\P(\F_{\S+}X_\T) \ge \F\P (X_\S).
$$
Since this holds for every $\P\in\mf{P}_{\S+},$ we have that 
$$
\F_{S+}X_\T \ge  X_\S.
$$
This proves the theorem if $\S, \T$ are optional times. In the case that they are stopping times, the theorem holds since the inequalies hold for all $\P\in\mf{P}_\S.$
\qed

%$\mathbbold{1}$

%J.J. Grobler, C.C.A. Labuschagne, Itô’s rule and Lévy’s theorem in vector lattices, J. %Math.
%Anal. Appl. (2017), http://dx.doi.org/10.1016/j.jmaa.2017.06.011


\begin{thebibliography}{99}
\bibitem{AB1} Aliprantis, C.D. and Burkinshaw, O., \textit{Locally solid Riesz spaces}, Academic Press, New York, San Francisco, London, 1978.
\bibitem{AB2} Aliprantis, C.D. and Burkinshaw, O., \textit{Positive Operators}, Academic Press Inc., Orlando, San Diego, New York, London, 1985.
%\bibitem{A} Ash, R.B., {\em Basic Probability Theory,} Wiley, New York, 1970.%
%\bibitem{AG} Ash, R.B. and Gardner, M.F., \textit{Topics in stochastic processes,} Academic Press Inc., New York, San Francisco London, 1975.%
%\bibitem{Ash} Ash, R.B. {\em Measure, Inegration and functional analysis}, Academic Press Inc., New York, 1972. 
%\bibitem{A&D} Ash, R.B. and Dol\'eans-Dade, C.A., {\em Probability \& Measure Theory}, Second Edition, Harcourt/Academic Press, Massachusetts, 2000.%

%\bibitem{B} Bartle, R.G. {A general bilinear vector integral,} \textit{Studia Math.}       %\textbf{15} (1956), 337-352. 
%\bibitem{AR1} Azouzi, Y. and Ramdane, K., {Burkholder Inequalities in Riesz spaces},%
%\textit{Indagationes Mathematicae,} \textbf{28} (2017), 1076-1094.%
%\bibitem{AR2}Azouzi, Y. and Ramdane, K., On the distribution function with respect to conditional expectation on Riesz spaces, \textit{Questiones Mathematicae}, doi/10.2989/16073606.2017.1377310.%
%\bibitem{BG} Blumenthal, R.M. and Getoor, R.D., {\em Markov processes and Potential Theory,} Academic Press Inc., New York, London, 1968. %
%\bibitem{CM} Cameron, R.H. and Martin, W.P., Transformation of Wiener integrals under translations, 
%\textit{Ann. Math.} \textbf{45} (1944), 386-396.%
%\bibitem{CL1} Cullender, S.F. and Labuschagne, C.C.A., {A description of norm-%convergent martingales on vector-valued $L^p$-spaces,} \textit{J. Math. Anal. Appl.}       %%\textbf{323} (2006), 119-130.
%\bibitem{CL2} Cullender, S.F. and Labuschagne, C.C.A., {Convergent martingales of %operators and the Radon-Nikod\`ym property in Banach spaces,} \textit{Proc. Amer. Math. %Soc.} \textbf{136} (2008), 3883-3893.
%\bibitem{DeMarr} DeMarr, R., A martingale convergence theorem in vector lattices,%textit{Canad. J. Math.} \textbf{18} (1966), 424--432.
\bibitem{DU} Diestel, J., J.J. Uhl, \textit{Vector measures}, Amer. Math. Soc., Providence, RI, 1977. 
%\bibitem{ID} Dobrakov, I., {On integration in Banach spaces,} \textit{Czech. Math, J.} \textbf{20}(3) (1970), %511-536.
\bibitem{PJD}
        P.J. Daniell,
        \emph{A general form of integral},
        Annals of Mathematics, Second Series, Vol. 19, no. 4 (Jun., 1918),
        279-294. (https://www.jstor.org/stable/1967495)
\bibitem{Donner} Donner, K., {\em Extension of Positive Operators and Korovkin Theorems}, Lecture Notes in Mathematics, Volume 904, Springer-Verlag, Berlin, Heidelberg, New York, 1982.
\bibitem{F} Fremlin, D.H., \textit{Topological Riesz spaces and measure theory,} Cambridge University Press, Cambridge, 1974.
%\bibitem{GTX} Gao, N., Troitsky, V.G., Xanthos, F., UO-convergence and its applications %to Ces\`aro means in Banach lattices, arXiv:1509.07914v1 [Math. FA] 25 Sep 2015 (To %appear in Israel J. Math).
%\bibitem{Gi} Girsanov, I.V., {On transforming a certain class of stochastic processes by absolutely continuous substitution of measures,} \textit{Theory Probab. Appl.} \textbf{5} (1960), 285-301.%
%\bibitem{G0} Grobler, J.J., {On the functional calculus in Archimedean Riesz spaces with applications to approximation theorems,} \textit{Quaestiones Math.} \textbf{11} (1988), 307-321.%
\bibitem{G1} Grobler, J.J., {Continuous stochastic processes in Riesz spaces: the Doob-Meyer decomposition}, \textit{Positivity} \textbf{14} (2010), 731-751.
\bibitem{G2} Grobler, J.J., {Doob's optional sampling theorem in Riesz spaces}, \textit{Positivity} \textbf{15} (2011), 617-637.
%\bibitem{G3} Grobler, J.J., {Jensen's and Martingale Inequalities in Riesz spaces}, \textit{Indagationes Mathematicae}, \textbf{25} (2014), 275-295 {dx.doi.org/10.1016/indag.2013.02.003}.%
%\bibitem{G4}  Grobler, J.J., {The Kolmogorov-\u{C}entsov theorem and Brownian motion in vector lattices}, \textit{J. Math. Anal. Appl.}, \textbf{410} (2014) 891-901.%
%\bibitem{G5} Grobler, J.J., Corrigendum to “The Kolmogorov-\u{C}entsov theorem and Brownian motion in vector lattices” [J. Math. Anal. Appl. 410 (2014) 891-901], DOI: 10.1016/j.jmaa.2014.05.068. %
\bibitem{G12} Grobler, J.J., {101 Years of vector lattice theory. A vector-valued Daniell integral} \textit{Preprint: Research Gate} 2020.
\bibitem{G6} Grobler, J.J., C.C.A. Labuschagne and Maraffa, V., {Quadratic variation of martingales in Riesz spaces}, \textit{J. Math. Anal. Appl.}, \textbf{410} (2014) 418-426.
\bibitem{G7} Grobler, J.J. and  C.C.A. Labuschagne, {The It\^o integral for Brownian motion in vector lattices: Part 1}, \textit{J. Math. Anal. Appl.}, \textbf{423} (2015) 797-819. doi 10.1016/j.jmaa.2014.08.013.
\bibitem{G8} Grobler, J.J. and  C.C.A. Labuschagne, {The It\^o integral for Brownian motion in vector lattices: Part 2}, \textit{J. Math. Anal. Appl.}, \textbf{423} (2015) 820-833. doi 10.1016/j.jmaa.2014.09.063.
%\bibitem{G10} Grobler, J.J. and  C.C.A. Labuschagne, {The quadratic variation of continuous time stochastic processes in vector lattices}, \textit{J. Math. Anal. Appl.}, \textbf{450} (2017) 314-329. doi 10.1016/j.jmaa.2017.01.034.%
\bibitem{G9} Grobler, J.J. and  C.C.A. Labuschagne, {The It\^o integral for martingales in vector lattices}, \textit{J. Math. Anal. Appl.} \textbf{450} (2017) 1245-1274. doi 10.1016/j.jmaa.2017.01.081.
%\bibitem{G11} Grobler, J.J. and C.C.A. Labuschagne, It\^o's rule and L\'evy's theorem in vector lattices, J. Math.%
%Anal. Appl. (2017), http://dx.doi.org/10.1016/j.jmaa.2017.06.011%
\bibitem{KS} Karatzas, I. and Shreve, S.E., \textit{Brownian motion and stochastic calculus}, Graduate Texts in Mathematics, Springer, New York, Berlin, Heidelberg, 1991.
%\bibitem{KW} Kunita, H. and Watanabe, S., {On square-integrable martingales}, \textit{Nagoya Math. J.} \textbf{30} (1967) 209--245. %
%\bibitem{Kuo} Kuo, H.-H. \textit{Introduction to Stochastic Integration,} Springer, New York, 2006.%
%\bibitem{Kuo} Kuo, W.-C, \textit{Stochastic processes on vector lattices,} Thesis, %University of the Witwatersrand, 2006.
%\bibitem{KLW1} Kuo, W.-C, Labuschagne, C.C.A., Watson, B.A., Discrete time stochastic %processes on Riesz spaces, \textit{Indag. Math.} \textbf{15} (2004), 435--451. 
%\bibitem{KLW2} Kuo, W.-C, Labuschagne, C.C.A., Watson, B.A.,  An upcrossing theorem for martingales on Riesz %spaces, \textit{Soft methodology and random information systems,} Springer-Verlag, 2004, 101--108. 
%\bibitem{KLW3} Kuo, W.-C, Labuschagne, C.C.A., Watson, B.A.,  Conditional expectation %on Riesz spaces, \textit{J. Math. Anal. Appl.,} \textbf{303} (2005), 509--521. 
%\bibitem{KLW4} Kuo, W.-C, Labuschagne, C.C.A., Watson, B.A.,  Zero-one laws for Riesz space and fuzzy random %variables, \textit{Fuzzy logic, soft computing and computational intelligence} Springer-Verlag and Tsinghua %University Press, Beijing, China, 2005, 393--397.   
%\bibitem{KLW5} Kuo, W.-C, Labuschagne, C.C.A., Watson, B.A., Convergence of Riesz space martingales, %\textit{Indag. Math.} \textbf{17} (2006), 271--283. 
%\bibitem{KLW6} Kuo, W.-C, Labuschagne, C.C.A., Watson, B.A., Ergodic theory and the %strong law of large numbers on Riesz spaces, \textit{J. Math. Anal. Appl.,}
%  \textbf{325} (2007), 422--437. 
%\bibitem{LW}  Labuschagne, C.C.A., Watson, B.A., Discrete time stochastic integrals in Riesz spaces, Positivity \textbf{14} (2010), 859-875.%
%\bibitem{L} L\'evy, P.  {\em Processus Stochastiques et Mouvement Brownien.}  Gauthier-Villars,  Paris 1948.%
%\bibitem{Loeve} Lo\'eve, M. \textit{Probabiblity Theory I and II} 4th Edition, Graduate Texts in Mathematics Volume 45, Springer-Verlag, New York, Berlin, Heidelberg, 1977.%
\bibitem{LZ} Luxemburg, W.A.J. and Zaanen, A.C., \textit{Riesz Spaces I}, North-Holland Publishing Company, Amsterdam, London, 1971.
%\bibitem{M} Meyer, P.A. {A decomposition theorem for supermartingales,} \textit{Illinois J. Math.} \textbf{6}%(1962), 193--205. 
\bibitem{MN} Meyer-Nieberg, P., \textit{Banach Lattices,} Springer-Verlag, Berlin, Heidelberg, New York, 1991.
\bibitem{Pr} Protter, P.E.,  \textit{Stochastic integration and Differential equations,} Sprinter-Verlag, Berlin, Heidelberg, New York, 2005.
%\bibitem{RY} Revuz, D. and Yor, M., \textit{Continuous Martingales and Brownian motion,} Springer-Verlag, Berlin, Heidelberg, New York, 1991.%
\bibitem{Sch} Schaefer, H.H., \textit{Banach lattices and positive operators,} Springer-Verlag, Berlin, Heidelberg, New York, 1974.
%\bibitem{SP} Schilling, R.L. and Partzsch, L., \textit{Brownian Motion, An Introduction to Stochastic Processes}, Walter de Gruyter GmbH \& Co, Berlin/Boston, 2012.%
%\bibitem{Sh} Shreve, S.E., \textit{Stochastic calculus for finance II, Continuous-time models,} Springer, New York, 2004.%
%\bibitem{Stoica1} Stoica, G., Vector-valued quasi-martingales, \textit{Stud. Cerc. %Mat.} \textbf{42} (1990), 73--79.
%\bibitem{Stoica2} Stoica, G., The structure of stochastic processes in normed vector %lattices, \textit{Stud. Cerc. Mat.} \textbf{46} (1994), 477--486.
%\bibitem{Troitsky} Troitsky, V., Martingales in Banach lattices, \textit{Positivity} \textbf{9} (2005), 437--456.%
%\bibitem{Va} Vardy, J.J., \textit{Markov processes and martingale generalisations on %Riesz spaces}, Ph.D.-thesis, University of the Witwatersrand, Johannesburg, South %Africa, 2012.%
%\bibitem{VW} Vardy, J.J. and Watson, B.A., Markov processes on Riesz spaces,  
%\textit{Positivity} \textbf{16} (2012), 373-391.%
\bibitem{Vu} Vulikh, B.Z., \textit{Introduction to the theory of partially ordered spaces,} Wolters-Noordhoff Scientific Publishers, Groningen, 1967.
\bibitem{W} Watson, B.A., An \^Ando-Douglas type theorem in Riesz spaces with a conditional expectation, {\em Positivity} \textbf{13} (2009), no 3, 543--558
\bibitem{Z1} Zaanen, A.C., \textit{Riesz spaces II}, North-Holland, Amsterdam, New York, 1983.
\bibitem{Z2} Zaanen, A.C., \textit{Introduction to Operator theory in Riesz spaces}, Springer-Verlag, Berlin, Heidelberg, New York, 1991.
%\bibitem{Z3} Zaanen, A.C., \textit{Continuity, Integration and Fourier theory}, Springer-Verlag, Berlin, Heidelberg, New York, London, Paris, Tokyo, 1989.% 
\end{thebibliography}
\end{document}